\documentclass{amsart}

\usepackage[numbers]{natbib}
\usepackage{rotating}
\usepackage{graphicx}

\newtheorem{theorem}{Theorem}
\newtheorem{example}[theorem]{Example}

\def\qspline #1{S\!(#1;q)} 
\def\iqspline #1#2{S_{\!#1}(#2;q)}

\def\D{D_q}
\def\hh #1{\mathcal{H}_ #1(q)}
\def\nn{\nonumber}

\begin{document}

\title{On Construction of A New Interpolation Tool: Cubic $q$-Spline}

\author{Orli Herscovici}

\address{O.~Herscovici\\Department of Mathematics,
University of Haifa,
3498838  Haifa, Israel\newline
Department of Mathematics, ORT Braude College, 2161002 Karmiel, Israel}
\email{orli.herscovici@gmail.com}

\maketitle

\begin{abstract}
This work presents a new interpolation tool, namely, cubic $q$-spline. Our new analogue generalizes a well known classical cubic spline. This analogue, based on the Jackson $q$-derivative, replaces an interpolating  piecewise cubic polynomial function by $q$-polynomials of degree three at most. The parameter $q$ provides a solution flexibility.
\end{abstract}

\section{INTRODUCTION}\vspace{-0.1cm}
The interpolation problem of an unknown function $f(x)$ 
when only the values of $f(x_i)$ at some point $x_i$ are 
given arises in different areas. One of 
widely used methods is a spline interpolation, and, 
particularly, cubic spline interpolation. That means that the 
function $f(x)$ is interpolated between two adjacent points 
$x_i$ and $x_{i+1}$ by a polynomial of degree three at most. Such 
interpolation is very suitable for smooth functions which do not have 
oscillating behaviour (cf. \cite{Atkinson1989}). Another advantage of the cubic 
polynomial interpolation is that it leads to a system of linear equations which is described by a tridiagonal matrix. This linear system has a unique solution which can be fast obtained. A generalization of cubic spline interpolation  was done by Marsden (see \cite{Marsden1982}). He chose a $q$-representation for knots and normalized $B$-splines as interpolating functions. $B$-splines are the generalization of the Bezier curves which are built with the help of Bernstein polynomials and they play an important role in the theory of polynomial interpolation. Their $q$-analogues were defined and studied in \cite{Budakci2015, Simeonov2013}. For further information and other details related to q-analogues of Bernstein polynomials, Bezier curves and splines, see the recent works \cite{Goldman2014, Goldman2015, Goldman2016, Oruc2003, Ostrovska2003, Simeonov2012}. Another $q$-generalization of polynomial interpolation was studied in \cite{Simsek2010}. Despite the popularity of $B$-spline, new generalizations of spline continue to appear. Some of them concerns about preserving of convexity \cite{Han2015}, some of them about smoothness of interpolation \cite{Segeth2018}, and others about degrees of freedom \cite{Mohammed2016}. Classical cubic spline already proved itself in data and curve fitting problems. Our new $q$-generalization gives it a new interesting twist. The matrix describing the linear system is tridiagonal like in the classic case, and a solution of this linear system can be obtained by simple recursive algorithm (see for example \cite{Kilic2008}). The $q$-parameter provides a flexibility of the solution.

We start from a short review of definitions coming from the 
quantum calculus (cf. \cite{Kac2002}).
The $q$-derivative is given by 
\[
\D f(x)=\frac{f(qx)-f(x)}{qx-x}. \nn
\]

\noindent For any complex number $c$, its $q$-analogue is  
defined as $
[c]_q=\frac{q^c-1}{q-1}$. For natural $n$, a $q$-factorial is defined as
$[n]_q!=\prod_{k=1}^{n}[k]_q$, with $ [0]_q!=1$.

The $q$-analogue of the polynomial is $(x-c)^n_q=\prod_{k=1}^n(x-cq^{k-1})$ by assuming as usually that $(x-c)^0_q=1$. 
It is easy to show that
$
\D (x-c)^n_q
=[n]_q(x-c)^{n-1}_q.
$
We will denote by $D^{k}_q$ the $k$-th $q$-derivative.
The Jackson $q$-integral of $f(x)$ is defined as (cf. \cite{Kac2002,Jackson1910})
\[
\int f(x)d_qx=(1-q)x\sum_{j=0}^{\infty}q^jf(q^jx), 
\label{q-integral}
\]
and $\D\int f(x)d_qx=f(x)$. Note that $q$ is usually considered to be $0<q<1$.

In the next section we describe a process of building of a cubic $q$-spline.

\section{BUILDING OF A $q$-ANALOGUE OF CUBIC SPLINE}
Let a function $f(x)$ is given by its values $f(x_i)=f_i$,  
$0\leq i \leq n$, at $n+1$ nodes $a=x_0<x_1<\ldots<x_n=b$.
We define the cubic $q$-spline function $\qspline x$ 
as a function of a variable $x$ with parameter $q$ as 
following: \vspace{-0.1cm}
\begin{equation}
\qspline x=\left\lbrace 
\begin{array}{rl}
\iqspline 1 x & x_0\leq x \leq x_1,\\
\iqspline 2 x & x_1< x \leq x_2,\\
\ldots & \ldots \\
\iqspline n x & x_{n-1}< x \leq x_n,
\end{array}\right.  \label{Sxq}
\end{equation}
where each $\iqspline i x$, $1\leq i \leq n$, is a $q$-polynomial in variable $x$
of degree at most three, and $\qspline x$, $\D\qspline x$, 
$\D^2\qspline x$ are 
continuous on $[a,b]$.
To provide these properties we demand 
\begin{equation}
\left\{\begin{array}{ll}
\iqspline i {x_{i-1}}=f(x_{i-1}), & i=1,\ldots,n,\\
\iqspline i {x_i}=f(x_i), & i=1,\ldots,n,
\end{array}\right.
\label{spline-cond}
\end{equation} 
\begin{eqnarray}
&\D \iqspline i{x_i}=\D \iqspline {{i+1}} {x_i}, 
& i=1,\ldots,n-1,\label{dif1-qspline-condition}\\
&\D^2 \iqspline i {x_i}=
\D^2 \iqspline {{i+1}}{x_i}, 
& i=1,\ldots,n-1.\label{dif2-qspline-condition}
\end{eqnarray}
The boundary conditions for the clamped cubic $q$-spline are
\begin{equation}
\left\{\begin{array}{ll}
\D \iqspline 1{x_0}=\D f(x_0),\\
\D \iqspline n{x_n}=\D f(x_n).
\end{array}\right.
\label{BCclamped}
\end{equation}

Let us denote by $\mu_i(q)$, $i=0,\ldots,n$ 
the value of the second
$q$-derivative of the spline $\qspline x$  at the node $x_i$,
that is $\mu_i(q)=\D^2\qspline {x_i}$. Since the spline 
$\qspline x$ is a polynomial of degree at most three, 
its second derivative is a polynomial of degree at most one. 
Therefore it can be written as
\[
\D^2\iqspline{i}{x}  
=\mu_{i-1}(q)\frac{x_i-x}{h_i}+\mu_i(q)\frac{x-x_{i-1}}{h_i},
\]
where $h_i=x_i-x_{i-1}$ and $x_{i-1}\leq x\leq x_i$ for 
$1\leq i\leq n$.
By performing $q$-integration we obtain
\begin{equation}	
\D\iqspline{i}{x}=\int \D^2\iqspline{i}{x}d_qx =
\frac{\mu_i(q)}{[2]_qh_i}(x-x_{i-1})^2_q
-\frac{\mu_{i-1}(q)}{[2]_qh_i}(x-x_i)^2_q+A_i(q).\label{FirstDer1}
\end{equation}
Let us denote by $f[x_{i-1},x_i]=\frac{f(x_i)-f(x_{i-1})}{x_i-x_{i-1}}$ the first order divided difference of  a function $f(x)$  at nodes $x_{i-1},x_i$. By integrating  (\ref{FirstDer1}), we obtain
\begin{equation}
\iqspline{i}{x}=\int \D\iqspline{i}{x}d_qx
=
\frac{\mu_i(q)}{[3]_q!h_i}(x-x_{i-1})^3_q
-\frac{\mu_{i-1}(q)}{[3]_q!h_i}(x-x_i)^3_q
+A_i(q)(x-x_{i-1})+B_i(q),\label{spline1}
\end{equation}
where $A_i(q)$ and $B_i(q)$ depend on parameter $q$ only and can be
found by substituting $x=x_{i}$ and $x=x_{i-1}$ in (\ref{spline1}) respectively
and using the conditions (\ref{spline-cond})  as following
\begin{eqnarray}
B_i(q) &=&f(x_{i-1})+\frac{\mu_{i-1}(q)}{[3]_q!h_i}(x_{i-1}-x_i)^3_q,
\label{Bi}\\
A_i(q)&=&\frac{f(x_i)-f(x_{i-1})}{h_i}
-\frac{\mu_i(q)}{[3]_q!h_i^2}(x_i-x_{i-1})^3_q
-\frac{\mu_{i-1}(q)}{[3]_q!h_i^2}(x_{i-1}-x_i)^3_q.
\label{Ai}
\end{eqnarray}
By substituting the detailed expressions (\ref{Bi}--\ref{Ai}) 
for functions $A_i(q)$ and $B_i(q)$, we obtain the $i$th spline function
as following
\begin{eqnarray}
\iqspline{i}{x}&=&\frac{\mu_i(q)}{[3]_q!h_i}(x-x_{i-1})^3_q
-\frac{\mu_{i-1}(q)}{[3]_q!h_i}(x-x_i)^3_q\nn\\
&+&\left(f[x_{i-1},x_i]-\frac{\mu_i(q)}{[3]_q!h_i^2}
(x_i-x_{i-1})^3_q-\frac{\mu_{i-1}(q)}{[3]_q!h^2_i}
(x_{i-1}-x_i)^3_q\right)(x-x_{i-1})\label{spline-formula}\\
&+ &f(x_{i-1})+\frac{\mu_{i-1}(q)}{[3]_q!h_i}(x_{i-1}-x_i)^3_q.\nn
\end{eqnarray}

In order to obtain the unknown moments $\mu_i(q)$, 
$i=0,\ldots,n$ 
we use the conditions (\ref{dif1-qspline-condition}) for the first 
derivative of the $q$-spline (\ref{FirstDer1}).
Hence, for $i=1,\dots,n-1$, we have
\begin{eqnarray}
\frac{\mu_i(q)}{[2]_qh_i}(x_i-x_{i-1})^2_q&+&f[x_{i-1},x_i]
-\frac{\mu_i(q)}{[3]_q!h_i^2}(x_i-x_{i-1})^3_q
-\frac{\mu_{i-1}(q)}{[3]_q!h_i^2}(x_{i-1}-x_i)^3_q\nn\\
&=&-\frac{\mu_i(q)}{[2]_qh_{i+1}}(x_i-x_{i+1})^2_q+f[x_i,x_{i+1}]
-\frac{\mu_{i+1}(q)}{[3]_q!h_{i+1}^2}(x_{i+1}-x_i)^3_q\nn\\
&-&\frac{\mu_i(q)}{[3]_q!h_{i+1}^2}(x_i-x_{i+1})^3_q.\label{mu_cond}
\end{eqnarray}

Moreover, from the boundary conditions (\ref{BCclamped}) for the clamped $q$-spline 
we obtain
\begin{eqnarray}
\D\iqspline{1}{x_0}&=&-\frac{\mu_0(q)}{[2]_qh_1}(x_0-x_1)^2_q
+f[x_0,x_1]\nn\\
&-&\frac{\mu_1(q)}{[3]_q!h_1^2}(x_1-x_0)^3_q-\frac{\mu_0(q)}{[3]_q!h_1^2}
(x_0-x_1)^3_q=\D f(x_0),\label{mu_cond0}\\
\D\iqspline{n}{x_n}&=&\frac{\mu_n(q)}{[2]_qh_n}(x_n-x_{n-1})^2_q
+f[x_{n-1},x_n]\nn\\&-&\frac{\mu_n(q)}{[3]_q!h^2_n}(x_n-x_{n-1})^3_q
-\frac{\mu_{n-1}(q)}{[3]_q!h_n^2}(x_{n-1}-x_n)^3_q=
\D f(x_n).\label{mu_condn}
\end{eqnarray}

The equations (\ref{mu_cond}--\ref{mu_condn}) form a system of 
$(n+1)$ linear equations with respect to the moments $\mu_i(q)$
that can be written shortly as $A\mu=b$.
With the notation $\hh{i}=\frac{(qx_i-x_{i-1})(x_i-qx_{i-1})}
{x_i-x_{i-1}}$, we have

\begin{equation}
A=
\left(
\begin{array}{cccccc}
[2]_q \hh{1} & \frac{(x_1-qx_0)^2_q}{x_1-x_0} & 0 & 0 & \ldots & 0 \\
\frac{(x_0-qx_1)^2}{x_1-x_0} & [2]_q(\hh{1}+\hh{2}) & 
\frac{(x_2-qx_1)^2_q}{x_2-x_1} & 0 & \ldots & 0\\
\ldots & \ldots & \ldots & \ldots & \ldots & \ldots\\
0 & 0 & 0 & 0 & \frac{(x_{n-1}-qx_n)^2_q}{x_n-x_{n-1}} & [2]_q\hh{n}
\end{array}\right),\label{matrix_A}
\end{equation}

\begin{equation}
b=[3]_q!\left(
\begin{array}{c}
f[x_0,x_1]-\D f(x_0) \\
f[x_1,x_2]-f[x_0,x_1]\\
\vdots\\
f[x_{n-1},x_n]-f[x_{n-2},x_{n-1}]\\
\D f(x_n)-f[x_{n-1},x_n]
\end{array}\right),\quad\quad
\mu=\left(
\begin{array}{c}
\mu_0(q) \\
\vdots\\
\mu_n(q)
\end{array}\right).
\label{vec_b}
\end{equation}

The unknown moments $\mu=\mu(q)$ (\ref{vec_b}) are the solution of the matrix 
equation $A\mu=b$, where $A$ is given by (\ref{matrix_A})
and $b$ is given by (\ref{vec_b}). In view of the fact that $A$ 
is a tridiagonal matrix, the method proposed in \cite{Kilic2008}
may be used for evaluating the moments $\mu_i(q)$, $i=0..,n$.

The moments $\mu_i(q)$, $i=0,\ldots,n$ are the solution of the 
matrix equation $A\mu=b$, therefore they are the unique functions 
of parameter $q$. 
Thus we can state the following result.
\begin{theorem}
Let the piecewise function $\qspline x$  be
defined by (\ref{Sxq}) and (\ref{spline-formula}), where $\mu=\mu(q)$ 
is a unique solution of the matrix equation $A\mu=b$, with $A$ given by (\ref{matrix_A}) and $b$ given by (\ref{vec_b}). Then  $\qspline x$ is a $q$-analogue of the cubic spline interpolation of a function $f(x)$.
\end{theorem}

\begin{example}
Let us consider a $q$-spline interpolation of a function $f(x)=x^4$ on the interval $x\in [-1,1]$ at the knots $x_0=-1$, $x_1=0$, $x_2=1$. We assume that the values of the function itself and its first derivation at the knots are given. The classical cubic spline solution is given in Example 1 on page 824--825 of \cite{Kreyszig2011}. By using the above proposed method,  one can obtain the cubic $q$-spline solution:
\[
S(x;q)=\left\{\begin{array}{lcl}\frac{q^3-3q^2-2q-2}{2q^2+q}x^3+\frac{q^3-q^2-q-2}{2q^2+q}x^2, &  &-1\leq x\leq 0,\\
-\frac{q^3-3q^2-2q-2}{2q^2+q}x^3+\frac{q^3-q^2-q-2}{2q^2+q}x^2, &  & 0\leq x\leq 1.\end{array}\right.
\]  
It is easy to obtain the classical cubic spline solution corresponding to $q\rightarrow1$. One can see that there exist  a slight oscillation of the spline regarding to the original function. Decreasing the parameter $q$ leads to more intensive oscillation. However, \textbf{increasing} the parameter $q$ overcomes the oscillation effect and interpolates the original function with significant improvement.

\begin{figure}[h]
\centerline{\includegraphics[height=120pt]{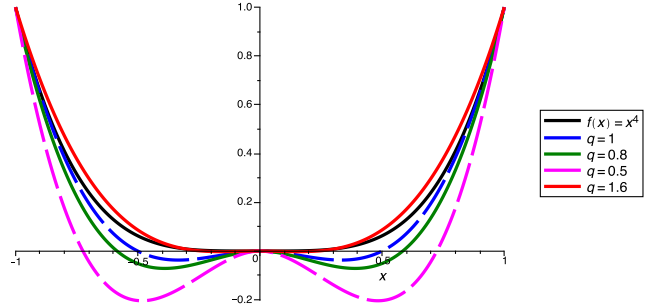}}
\caption{Interpolation of $f(x)=x^4$ by cubic $q$-spline.}
\end{figure}
\end{example}

\section{ACKNOWLEDGMENTS}
This research 
was supported  by the Ministry of Science and Technology,
Israel.

\noindent
The author thanks to the anonymous referees for their advices.

\nocite{*}
\bibliographystyle{aipnum-cp}
\bibliography{sample}

\end{document}